\documentclass[12pt]{article}
\usepackage[cp1251]{inputenc}
\usepackage{amsfonts,amsmath}
\newtheorem{theorem}{Theorem}
\newtheorem{cons}{Corollary}
\newtheorem{rema}{Remark}

\begin{document}

\title {\hbox{\normalsize }\hbox{}
Positive one--point commuting difference operators \thanks{The research was supported by the Russian Foundation for Basic Research (Grant 18-01-00411)}}

\author {Gulnara S.~Mauleshova, \ Andrey E.~Mironov}

\date{}
\maketitle


\sloppy
\begin{abstract}
In this paper we study commuting difference operators containing a shift operator with only positive degrees. We construct examples of such operators in the case of hyperelliptic spectral curves.
\end{abstract}

\section{Introduction and Main Results}

Every maximal commutative ring of difference operators is isomorphic to a ring of meromorphic functions on an algebraic spectral curve $\Gamma$ with poles $p_1, \ldots, p_s.$ Such operators are called {\it $s$--point operators}. Common eigenfunctions of these operators form a vector bundle of rank $l$ over an affine part of $\Gamma,$ where $l$ is called the {\it rank} of the operators. Krichever~\cite{K3} (see also~\cite{Mum}) found two--point operators of rank one. Krichever and Novikov~\cite{KN2},~\cite{KN3} essentially developed a classification of commuting rings of difference operators, but the classification is not completed. The existence of one--point higher rank operators is proved in~\cite{KN2},~\cite{KN3}. Such operators were studied~\cite{MM1} in the case of hyperelliptic spectral curves.  We mention here that in every maximal ring of operators discovered in~\cite{K3}--\cite{MM1},  there are operators of the form
$$
v_m(n)T^m+\ldots+v_0(n)+v_{-1}(n)T^{-1}+\ldots+v_{-k}(n)T^{-k}, \quad m,k>0,
$$
where $T\psi(n)=\psi(n+1),$ i.e. the shift operator is included in the operator with positive as well as negative degrees.
In our short note~\cite{MM2} we found a new class of rank one commuting difference operators such that these operators have the form
$$
L_m = T^m+u_{m-1}(n)T^{m-1}+\ldots+u_0(n),
$$
i.e. $T$ has only positive degrees in $L_m.$ We call such operators {\it positive operators}.

In this paper we give proofs of the results announced in~\cite{MM2}, construct new examples (Theorem 5, see below), and also generalise the construction from~\cite{MM2} to the higher rank operators.

We choose the following spectral data
$$
S=\{\Gamma, \gamma_1,\ldots,\gamma_g, q, k^{-1}, K_n\},
$$
where $\Gamma$ is a Riemann surface of genus $g$, $\gamma = \gamma_1 + \cdots + \gamma_g$ is a non--special divisor on $\Gamma,$ $q\in\Gamma $ is a fixed point, $k^{-1}$ is a local parameter near $q$, $K_n\in\Gamma$, $n\in\mathbb Z$ is a set of points in general position. If $K_n$ does not depend on $n$, i.e. $K_n=q^+\neq q$ then we get the spectral data for two--point operators of rank one (see~\cite{K3}).
\begin{theorem}[\cite{MM2}]
There exists a unique function $\psi(n,P), \ n\in \mathbb{Z}, \ P \in \Gamma$, which has the following properties:

\noindent{\bf 1}. The divisor of zeros and poles of $\psi$ has the form
$$
\gamma_1(n)+\ldots+\gamma_g(n)+K_1+\ldots+K_n-\gamma_1-\ldots-\gamma_g-nq,
$$
if $n > 0$ and has the form
$$
\gamma_1(n)+\ldots+\gamma_g(n)-K_{-1}-\ldots-K_{n}-\gamma_1-\ldots-\gamma_g-nq,
$$
if $n < 0.$

\noindent{\bf 2}. In the neighborhood of $q$, the function $\psi$ has the form
$$
\psi=k^n+O(k^{n-1}).
$$
\noindent{\bf 3}. $\psi(0,P) = 1.$
\end{theorem}
The function $\psi(n,P)$ is called the {\it Baker--Akhiezer function}.
For a meromorphic function $g(P)$ on $\Gamma$ with a unique pole of order $m$ in $q,$ $g = k^m + O(k^{m-1})$ there exists a unique operator such that $L_m \psi(n,P) = g(P) \psi(n,P)$.
\begin{rema}
In the case of the two--dimensional discrete Schr\"{o}dinger operator, Krichever~\cite{K5} considered spectral data in which an additional set of points $K_n$ appears similar to our construction.
\end{rema}

We consider the hyperelliptic spectral curve $\Gamma$ given by the equation
$$
w^2 = F_g(z) = z^{2g+1}+c_{2g} z^{2g}+\ldots+c_0,
$$
we choose $q=\infty.$ Let $\psi(n,P)$ be the corresponding Baker--Akhiezer function. There exist commuting operators $L_2,L_{2g+1}$ such that
$$
L_2\psi=((T+U_n)^2+W_n)\psi=z\psi, \quad L_{2g+1}\psi=w\psi.
$$
\begin{theorem}[\cite{MM2}]
The following identity holds:
$$
L_2-z=(T+U_n+U_{n+1}+\chi_n(P))(T-\chi_n(P)),
$$
where
$$
\chi=\frac{\psi(n+1,P)}{\psi(n,P)}=\frac{S_n}{Q_n}+\frac{w}{Q_n},
\quad
S_n(z)=-U_nz^g+\delta_{g-1}(n)z^{g-1}+\ldots+\delta_0(n),
$$
\begin{equation}
\label{eq1}
Q_n=-\frac{S_{n-1}+S_n}{U_{n-1}+U_n},
\end{equation}
Functions $U_n, W_n, S_n$ satisfy the equation
\begin{equation}
\label{eq2}
F_g(z)=S^2_n + (z-U^2_n-W_n)Q_{n}Q_{n+1}.
\end{equation}
\end{theorem}
The equation (\ref{eq2}) can be linearized. Namely, if we replace $n$ by $n + 1$ in (\ref{eq2}), take a difference with (\ref{eq2}) and apply (\ref{eq1}), then the right hand side is divisible by $Q_{n + 1}.$ So, we obtain the linear equation on $S_n.$

\begin{cons}[\cite{MM2}]
Functions $S_n(z), U_n, W_n$ satisfy the equation
\begin{equation}
\label{eq3}
(U_n+U_{n+1})(S_n-S_{n+1}) - (z-U^2_n-W_n)Q_n + (z-U^2_{n+1}-W_{n+1})Q_{n+2}=0.
\end{equation}
\end{cons}
In the case of the elliptic spectral curve $\Gamma$ given by the equation
$$
w^2=F_1(z)=z^3 + c_2 z^2 + c_1 z + c_0,
$$
the equation (\ref{eq2}) is exactly solvable.
\begin{cons}[\cite{MM2}]
Operator
$$
L_2=(T+U_n)^2+W_n,
$$
where
\begin{equation}
\label{eq4}
U_n=-\frac{\sqrt{F_1(\gamma_n)}+\sqrt{F_1(\gamma_{n+1})}}{\gamma_n-\gamma_{n+1}},\quad W_n=-c_2-\gamma_n-\gamma_{n+1},
\end{equation}
$\gamma_n$ is an arbitrary functional parameter, commutes with the operator
$$
L_3 = T^3 + \big(U_n + U_{n+1} + U_{n+2} \big) T^2 + \big(U^2_n + U^2_{n+1} + U_n U_{n+1} + W_n -\gamma_{n+2}\big) T +
$$
$$
+\big(\sqrt{F_1(\gamma_n)} + U_n (U^2_n + W_n - \gamma_n)\big).
$$
\end{cons}
Theorem 2 allows us to construct explicit examples.
\begin{theorem}[\cite{MM2}]
Operator
$$
L_2=(T+r_1\cos(n))^2+\frac{r_1^2 \sin(g)\sin(g+1)}{2 \cos^2(g+\frac{1}{2})}\cos(2n), \quad r_1\neq0
$$
commutes with an operator $L_{2g + 1}$ of order $2g+1$.
\end{theorem}

\begin{theorem}[\cite{MM2}]
Operator
$$
L_2=(T+\alpha_2n^2+\alpha_0)^2-g(g+1)\alpha_2^2n^2, \quad \alpha_2\neq0
$$
commutes with an operator $L_{2g + 1}$ of order $2g+1$.
\end{theorem}

\begin{theorem}
Operator
$$
L_2=(T+\beta_1 a^n)^2 + \frac{a^{2 g} + a^{2 g+2} - a^{4 g+2} - 1}{(a^{2 g+1} + 1)^2} \beta^2_1 a^{2 n}, \quad \beta_1\neq0, \quad a\neq1,
$$
commutes with an operator $L_{2g + 1}$ of order $2g+1$.
\end{theorem}

\begin{rema}
Theorem 2 has the following remarkable application.

Let $\wp(x), \ \zeta(x)$ are Weierstrass functions
$$
(\wp'(x))^2 = 4 \wp^3(x)+g_2 \wp(x)+g_3, \qquad -\zeta'(x) = \wp(x).
$$
Set
$$
A_1=-2 \zeta(\varepsilon)-\zeta(x-\varepsilon)+\zeta(x+\varepsilon),
$$
$$
A_2=-\frac{3}{2} (\zeta(\varepsilon)+\zeta(3 \varepsilon)+\zeta(x-2 \varepsilon)-\zeta(x+2 \varepsilon),
$$
Further, for odd $g=2 g_1+1$ we set
$$
A_g=A_1\prod_{k=1}^{g_1}\bigg(1+\frac{\zeta(x-(2 k+1) \varepsilon)-\zeta(x+(2 k+1) \varepsilon)}{\zeta(\varepsilon)+\zeta((4 k+1) \varepsilon)} \bigg),
$$
for even $g=2 g_1$ we set
$$
A_g=A_2\prod_{k=2}^{g_1}\bigg(1+\frac{\zeta(x-2 k \varepsilon)-\zeta(x+2 k \varepsilon)}{\zeta(\varepsilon)+\zeta((4 k-1) \varepsilon)} \bigg).
$$
Set
$$
L_2 = \frac{T^2_\varepsilon}{\varepsilon^2}+A_g(x,\varepsilon) \frac{T_\varepsilon}{\varepsilon}+\wp(\varepsilon),
$$
where $T_\varepsilon$ is the shift operator, $T_\varepsilon \psi(x)= \psi(x+\varepsilon).$

\begin{theorem}[\cite{MM3}]
The operator $L_2$ commutes with an operator of the form
$$
L_{2g+1}=\sum_{j=0}^{2g+1} B_j(x,\varepsilon) \frac{T_\varepsilon^j}{\varepsilon^j}.
$$
There is a decomposition
$$
L_2 = \partial_x^2-g(g+1) \wp(x)+O(\varepsilon).
$$
\end{theorem}
We checked that when $g=1,2,$ and numerically checked that when $g=3,\ldots, 6$ the spectral curve of $L_2, \ L_{2g+1}$ does not depend on the parameter $ \varepsilon $ and coincides with the spectral curve of the Lam\'{e} operator $\partial_x^2-g(g+1) \wp(x)$. So, we have a remarkable discretization of the Lam\'{e} operator preserving all integrable properties. The proof of Theorem 6 is very long, so we do not present it in this paper.
\end{rema}
In Section 2 we give proofs of Theorems 1--5 and examples. In Section 3 we point out the spectral data for one--point positive commuting difference operators of higher rank.

\section{Proofs of Theorems 1--5}

In this section we give proofs of Theorems 1--5.

\subsection{Theorem 1}

Consider the case $n> 0$ (the case $n <0$ is similar). By the Riemann--Roch theorem, the dimension of the space of meromorphic functions on $\Gamma$ with the pole divisor $\gamma_1 + \ldots + \gamma_g + nq$ is equal to
$$
l(\gamma_1 + \ldots + \gamma_g + nq) = n+1.
$$
The requirement that a function from this space vanishes at points $K_1, \ldots, K_n$ extracts a one--dimensional subspace of this space. Condition (2) gives us the unique function. Theorem 1 is proved.

For any meromorphic function $g(P)$ on $\Gamma$ with the unique pole of order $m$ at $q$ there exists a unique operator of the form
$$
L_m = T^m + u_{m-1}(n)T^{m-1} + \ldots + u_0(n)
$$
such that
$$
L_m \psi = g(P) \psi.
$$
Indeed, consider the case $n> 0$ (the case $n <0$ is similar). Let us note that for any $u_j(n)$ the pole divisor of the function
$$
\varphi(n,P) = L_m \psi - g(P) \psi
$$
has the form
$$
\gamma_1 + \ldots + \gamma_g + (m+n-1) q.
$$
Moreover, the function $\varphi(n,P)$ has zeros at points $K_1, \ldots, K_n$. Therefore, we can choose functions $u_{m-1}(n), \ldots, u_0(n)$ such that the function $\varphi(n,P)$ has zeros at points $K_1, \ldots, K_n, K_{n+1}, \ldots, K_{m+n}.$ Thus $\varphi(n,P) = 0.$ Similarly, for a meromorphic function $f(P)$ with the unique pole of order $s$ at $q$ there exists an operator $L_s$ such that $L_s \psi = f(P) \psi.$ Operators $L_m$ and $L_s$ commute.

\subsection{Theorem 2}

We have the identity
$$
L_2 \psi = \big((T+U_n)^2+W_n\big) \psi = z \psi.
$$
Set
$
\chi_n(P)=\frac{\psi(n+1,P)}{\psi(n,P)},
$
then $(T-\chi_n(P)) \psi(n,P) = 0.$

The operator $L_2-z$ can be factorized
\begin{equation}
\label{eq5}
(T+U_n)^2+W_n-z=(T+\tilde{\chi}_{n})(T-\chi_n(P)),
\end{equation}
where
$$
\tilde{\chi}_{n} = U_n+U_{n+1}+\chi_{n+1}(P),
$$
herewith $\chi_n(P)$ satisfies the equation
\begin{equation}
\label{eq6}
-z + U_n^2 + W_n + \chi_n(P) (U_n + U_{n+1} + \chi_{n+1}(P))=0.
\end{equation}
The equation (\ref{eq6}) follows from the direct comparison of the left and right parts in (\ref{eq5}).
The function $\psi(n,P)$ satisfies (1)--(3), consequently, the function $\chi_n(P)=\frac{\psi(n+1,P)}{\psi(n,P)}$ has the first order pole in $q,$ as well as simple poles in $\gamma_1(n),\ldots,\gamma_g(n).$ Thus $\chi_n(P)$ is a rational function with the pole divisor
$$
\big(\chi_n(P)\big)_\infty = \gamma_1(n) + \ldots + \gamma_g(n) + q,
$$
and in the neighborhood of $q$ it has the expansion
$$
\chi_n(P) = \frac{1}{k} + O(1).
$$
Let $\gamma_i(n)$ have coordinates $(\alpha_i(n),\mu_i(n)).$
We denote the polynomial of degree $g$ with respect to $z$ by $Q_n$
$$
Q_n=(z-\alpha_1(n))\ldots(z-\alpha_g(n)).
$$
Then the rational function $\chi_n(P)$ has the form
\begin{equation}
\label{eq7}
\chi_n(P) = \frac{S_n}{Q_n} + \frac{w}{Q_n},
\end{equation}
where $S_n$ is a polynomial of degree $g$ with respect to $z$
$$
S_n=\delta_g(n) z^g+\delta_{g-1}(n)z^{g-1}+\ldots+\delta_0(n).
$$
Substituting (\ref{eq7}) in (\ref{eq6}), we get that the polynomial $Q_n$ has the form
$$
Q_n=-\frac{S_{n-1}+S_n}{U_{n-1}+U_n}.
$$
Hence $\delta_g(n) = -U_n.$
Thus
\begin{equation}
\label{eq8}
S_n=-U_n z^g+\delta_{g-1}(n)z^{g-1}+\ldots+\delta_0(n)
\end{equation}
satisfies the equation (\ref{eq2}). Theorem 2 is proved.

In the proof of theorems 3 and 4 we use the following observation. Let us assume that $U_n, \ W_n$ and $S_n$ are even functions in $n$
$$
U_n = U_{-n}, \quad W_n = W_{-n}, \quad S_n = S_{-n}.
$$
Let $R_n$ be the left hand side of (\ref{eq3})
$$
R_n = S_{n-1} (U_{n+1} + U_{n+2}) (z - U^2_n - W_n) +
$$
$$
 S_n (U_{n+1} + U_{n+2}) (z + U_n U_{n+1} + U_{n-1} (U_n + U_{n+1}) - W_n) -
$$
$$
 S_{n+1} (U_{n-1} + U_n) (z +U_n U_{n+1} + (U_n + U_{n+1}) U_{n+2} - W_{n+1}) -
$$
$$
 S_{n+2} (U_{n-1} + U_n) (z - U^2_{n+1} - W_{n+1}).
$$
Then $R_n$ is skew invariant under the replacement $n \rightarrow -n-1,$ i.e. $R_n=-R_{-n-1}.$

\subsection{Theorem 3}

To prove Theorem 3 it is enough to show that for
$$
U_n=r_1\cos(n), \qquad W_n=\frac{r_1^2\sin(g)\sin(g+1)}{2 \cos^2(g+\frac{1}{2})} \cos(2n)
$$
there is a polynomial (\ref{eq8}) satisfying (\ref{eq3}).
\\
In our case the equation (\ref{eq3}) has the form
$$
R_n = \frac{1}{2} r_1 S_{n-1} (\cos(n+1) + \cos(n+2)) \big(2 z - r_1^2 (2 \cos^2(n)+
$$
$$
\frac{\sin(g) \sin(g+1)}{\cos^2(g+\frac{1}{2})}  \cos(2 n))\big) +\frac{1}{2} r_1 S_n (\cos(n+1) + \cos(n+2)) \big(2 z +
$$
$$
r_1^2 (2 \cos(1) + \cos(2) +\frac{\cos^2(\frac{1}{2}) (2 \cos(2 g+1)+1)}{\cos^2(g+\frac{1}{2})} \cos(2 n))\big) -
$$
$$
r_1 S_{n+1} (\cos(1 - n) + \cos(n)) \big(z + \frac{1}{2} r_1^2 (\cos(1) + 2 \cos(n) (\cos(n+1) +
$$
$$
\cos(n+2)) + \cos(2 n+3) - \frac{\sin(g) \sin(g+1)}{\cos^2(g+\frac{1}{2})} \cos(2 (n+1)) \big) -
$$
$$
\frac{1}{2} r_1 S_{n+2} (\cos(1 - n) + \cos(n))\big(2 z - r_1^2 (2 \cos^2(n+1) +
$$
\begin{equation}
\label{eq9}
\frac{\sin(g) \sin(g+1)}{\cos^2(g+\frac{1}{2})} \cos(2 (n+1))) \big)=0.
\end{equation}
Let
$$
S_n=A_{2g+1}\cos((2g+1)n)+A_{2g-1}\cos((2g-1)n)+\ldots+A_3\cos(3n)+A_1\cos(n),
$$
where $A_i=A_i(z).$ Then $U_n, W_n$ and $S_n$ are even functions in $n,$ and, as it was mentioned above, $R_n=-R_{-n-1}.$ It is a remarkable fact that after substituting $S_n$ in (\ref{eq9}) we obtain
$$
R_n=\alpha_{2g+4}\cos((2g+4)n)+\alpha_{2g+2}\cos((2g+2)n)+\ldots+\alpha_{4}\cos(4n)+\alpha_{2}\cos(2n)=0,
$$
where
$$
\alpha_s=8r_1^2A_{s-3}\cos^2(\frac{1}{2})(\cos(2g+1)-\cos(3-s))\sin(\frac{3-s}{2})+
$$
$$
A_{s-1}\bigg(8z\sin(\frac{5-3s}{2})-r_1^2\sin(\frac{1-3s}{2})-4r_1^2\sin(\frac{5-3s}{2})+
$$
$$
4r_1^2\sin(\frac{1-s}{2})+8z\sin(\frac{3-s}{2})+4r_1^2\sin(\frac{5-s}{2})+2r_1^2\sin(\frac{7-s}{2})+
$$
$$
2\sin(1)\sin(2g)(-2(2z+r_1^2(-1+\cos(1)+\cos(2)))\cos(\frac{s}{2})\sin(\frac{3}{2})+
$$
$$
2(r_1^2-2z))\cos(\frac{3s}{2})\sin(\frac{5}{2})+2\cos(\frac{3}{2})(2z+r_1^2(1+3\cos(1)+
$$
$$
\cos(2)))\sin(\frac{s}{2})-2(r_1^2-2z))\cos(\frac{5}{2})\sin(\frac{3s}{2}))+
$$
$$
2\cos(1)\cos(2g)(2(2z+r_1^2(-1+\cos(1)+\cos(2)))\cos(\frac{s}{2})\sin(\frac{3}{2})-
$$
$$
2(r_1^2-2z)\cos(\frac{3s}{2})\sin(\frac{5}{2})-2\cos(\frac{3}{2})(2z+r_1^2(1+3\cos(1)+
$$
$$
\cos(2)))\sin(\frac{s}{2})+2(r_1^2-2z)\cos(\frac{5}{2})\sin(\frac{3s}{2}))-3r_1^2\sin(\frac{1+s}{2})+
$$
$$
r_1^2\sin(\frac{3(1+s)}{2})-2r_1^2\sin(\frac{3+s}{2})-r_1^2\sin(\frac{5+s}{2})+2r_1^2\sin(\frac{1+3s}{2})\bigg)-
$$
$$
A_{s+1}\bigg(2r_1^2\sin(\frac{1-3s}{2})+r_1^2\sin(\frac{3}{2}-\frac{3s}{2})-3r_1^2\sin(\frac{1-s}{2})-2r_1^2\sin(\frac{3-s}{2})-
$$
$$
r_1^2\sin(\frac{5-s}{2})+2\cos(1)\cos(2g)(2(2z+r_1^2(-1+\cos(1)+\cos(2)))\cos(\frac{s}{2})\sin(\frac{3}{2})-
$$
$$
2(r_1^2-2z)\cos(\frac{3s}{2})\sin(\frac{5}{2})+2\cos(\frac{3}{2})(2z+r_1^2(1+3\cos(1)+
$$
$$
\cos(2)))\sin(\frac{s}{2})-2(r_1^2-2z)\cos(\frac{5}{2})\sin(\frac{3s}{2}))+2\sin(1)\sin(2g)(-2(2z+
$$
$$
r_1^2(-1+\cos(1)+\cos(2)))\cos(\frac{s}{2})\sin(\frac{3}{2})+2(r_1^2-
$$
$$
2z)\cos(\frac{3s}{2})\sin(\frac{5}{2}-2\cos(\frac{3}{2})(2z+r_1^2(1+3\cos(1)+\cos(2)))\sin(\frac{s}{2})+
$$
$$
2(r_1^2-2z)\cos(\frac{5}{2})\sin(\frac{3s}{2})+4r_1^2\sin(\frac{1+s}{2})+8z\sin(\frac{3+s}{2})+
$$
$$
4r_1^2\sin(\frac{5+s}{2})+2r_1^2\sin(\frac{7+s}{2})-r_1^2\sin(\frac{1+3s}{2})-4r_1^2\sin(\frac{5+3s}{2})+
$$
$$
8z\sin(\frac{5+3s}{2}))+8r_1^2A_{s+3}\cos^2(\frac{1}{2})(-\cos(1+2g)+\cos(3+s))\sin(\frac{3+s}{2}).
$$
From the above formula it follows that
$$
\alpha_{2g+4}=-16 A_{2 g+7} r_1^2 \cos^2(\frac{1}{2})\sin(3) \sin(2g+4) \sin(g+\frac{7}{2}) -
$$
$$
A_{2 g+5} \big(r_1^2 \sin(\frac{1}{2} - g) + 2 r_1^2 \sin(\frac{3}{2} - g) + 16 z \cos^2(g+\frac{1}{2}) (\sin(g+\frac{7}{2}) + \sin(3g+\frac{17}{2}))+
$$
$$
r_1^2 (2 \sin(\frac{7}{2} - g) + \sin(\frac{9}{2}- g) + 2 \sin(g+\frac{3}{2}) +
    4 \sin(g+\frac{5}{2}) + 4 \sin(g+\frac{9}{2}) + 2 \sin(g+\frac{11}{2}) -
$$
$$
2 \sin(g+\frac{15}{2}) + \sin(3g+\frac{1}{2}) + 2 \sin(3g+\frac{3}{2}) +
2 \sin(3g+\frac{5}{2}) + 2 \sin(3g+\frac{7}{2}) -
$$
$$
\sin(3g+\frac{9}{2}) - 4 \sin(3g+\frac{17}{2}) - 2 \sin(5g+\frac{19}{2}))\big) -
$$
$$
2 A_{2g+3} \sin(2 g+2)  \big(16 z \cos^2(g+\frac{1}{2}) \cos(g+\frac{3}{2}) +
r_1^2 (2 \cos(\frac{3}{2} - g) + \cos(\frac{5}{2} - g) +
$$
$$
2 \cos(g+\frac{1}{2}) - 2 \cos(g+\frac{3}{2}) + 2 \cos(g+\frac{5}{2}) - 2 \cos(g+\frac{9}{2}) -
\cos(g+\frac{11}{2}) - 2 \cos(3g+\frac{5}{2})) \big)=0
$$
(we assume that $A_{2g+3}=A_{2g+5}=A_{2g+7}=0$). So, we have $g+1$ equations $\alpha_{2g+2}=\alpha_{2g}=\ldots=\alpha_{2}=0$ and $g+1$ unknown functions $A_1, A_3,\ldots, A_{2g+1}.$ From the equations
$$
\alpha_{2g+2}=\alpha_{2g}=\ldots=\alpha_{4}=0
$$
we express $A_{2g-1}, A_{2g-3},\ldots, A_3, A_1$ via $A_{2g+1}.$ For example,
$$
A_{2g-1} = - A_{2g+1}\frac{\cos(\frac{1}{2} + g) (4 z + r_1^2 (2 \cos(1) + \cos(2)-1) -
   2 (r_1^2 - 2 z) \cos(2 g+1))}{8 r_1^2 \cos^3(\frac{1}{2} \sin(\frac{1}{2}) \sin(\frac{1}{2}- g)},
$$
$$
A_{2g-3} = -\frac{\cos(g+\frac{1}{2})}{16 r_1^2 \cos^3(\frac{1}{2}) \sin(\frac{3}{2} - g) \sin(1 - 2 g) \sin(-\frac{1}{2}) \cos(1)} \times
$$
$$
\bigg(A_{2 g-1} \big((r_1^2 - 2 z) \cos(4 g) \sin(2)-(2 z + r_1^2 (2 \cos(1) + \cos(2))) \sin(2) +
$$
$$
 2 (r_1^2 + 2 z \cos(1)) \cos(2) \sin(2 g) +
2 (r_1^2 - 2 z) \cos(2 g) (2 \cos^2(1) \sin(1) - \cos(2) \sin(2 g)) \big) +
$$
$$
A_{2 g+1} \sin(2g+2) \big(4 z + r_1^2 (2 \cos(1) + \cos(2)-1) - 2 (r_1^2 - 2 z) \cos(2 g+1)\big)\bigg).
$$
It turns out that $\alpha_{2}=0$ automatically. Indeed, from
$$
R_n+R_{-n-1} = \alpha_{2g+4}(1+\cos(2g+4))\cos((2g+4)n)+\ldots+\alpha_{4}(1+\cos(4))\cos(4n)+
$$
$$
\alpha_{2}(1+\cos(2))\cos(2n) - \alpha_{2g+4} \sin(2g+4) \sin((2g+4)n)-\ldots-\alpha_{4} \sin(4) \sin(4n)-
$$
$$
\alpha_{2} \sin(2) \sin(2n) =\alpha_{2}(1+\cos(2))\cos(2n) -\alpha_{2} \sin(2) \sin(2n) = 0.
$$
Hence $\alpha_{2}=0.$
Next, we choose $A_{2g+1}$ such that coefficient at $z^g$ in $S_n(z)$ is $-U_n.$
We find $S_n,$ satisfying the equation (\ref{eq3}). Theorem 3 is proved.

\noindent{\bf Example 2}. At $g=1$ we have
$$
U_n=r_1\cos(n), \qquad W_n=\frac{r_1^2\sin(1)\sin(2)}{2 \cos^2(\frac{3}{2})} \cos(2n).
$$
$$
S_n = A_3(z) \cos(3 n) + A_1(z) \cos(n),
$$
$$
A_3(z) = \frac{r_1^3 \sin^2(\frac{1}{2})}{(2 \cos(1)-1)^3},
$$
$$
A_1(z) =-\frac{r_1 (2 z (1 - 2 \cos(1))^2 + r_1^2 (5 \cos(1) - 2 \cos(2)-3))}{2 (1 -2 \cos(1))^2}.
$$
$$
Q_n =z - \frac{2 r_1^2 (\cos(1 - 2 n)-2 \cos(1)) \sin^2(\frac{1}{2})}{(1 - 2 \cos(1))^2}.
$$
The spectral curve is given by the equation
$$
w^2 =\big(z - \frac{4 r_1^2 \sin^4(\frac{1}{2})}{(1 - 2 \cos(1))^2} \big)^2 \big(z - \frac{r_1^2 (1 - 2 \cos(1) + \cos(2))}{(1 - 2 \cos(1))^2} \big).
$$


\subsection{Theorem 4}

Let
$$
U_n=a_2n^2+a_0, \qquad W_n=-g(g+1)a_2^2n^2.
$$
In this case the equation (\ref{eq3}) has the form
$$
R_n = \big(S_{n - 1} + S_n\big) \big(2 a_0 + a_2 (5 + 6 n + 2 n^2)\big) \big(a_2^2 g (1 + g) n^2 + z -
$$
$$
(a_0 + a_2 n^2)^2\big) +\big(S_n - S_{n + 1}\big) \big(2 a_0 + a_2 - 2 a_2 n + 2 a_2 n^2\big) \big(2 a_0 + a_2 + $$
$$
2 a_2 n + 2 a_2 n^2\big)\big(2 a_0 + a_2 (5 + 6 n + 2 n^2)\big) -\big(S_{n + 1} + S_{n + 2}\big) \big(2 a_0 + a_2 - 2 a_2 n +
$$
\begin{equation}
\label{eq10}
2 a_2 n^2\big) \big(a_2^2 g (1 + g) (1 + n)^2 - (a_0 + a_2 (1 + n)^2)^2 + z\big)=0.
\end{equation}
Let
$$
S_n=B_{2g+2}n^{2g+2}+B_{2g}n^{2g}+\ldots+B_2n^2+B_0, \qquad B_i=B_i(z).
$$
Then $U_n, W_n$ and $S_n$ are even functions in $n,$ and, as it was mentioned above, $R_n=-R_{-n-1}.$
After substituting $S_n$ in (\ref{eq10}) we obtain
$$
R_n = \vartheta_{2g+8}(z)n^{2g+8}+\vartheta_{2g+7}(z)n^{2g+7}+\ldots+\vartheta_0(z)=0,
$$
where
$$
\vartheta_s=-2 B_{s-3}a_2^3(5+2g-s)(s-4)(2g+s-3)+
$$
$$
B_{s-2}a_2^3((s-4)(s-2)-4g(g+1))(s-3)(s+1)+
$$
$$
B_{s-1}\big(\frac{1}{6}a_2^2(a_2(3(s-3)(8+(s-3)s(8+(s-3)s))-4g(s-1)(27+s(5s-22))-
$$
$$
4g^2(s-1)(27+s(5s-22)))-12a_0(6+s(-29+4g(g+1)-3(s-6)s)))-8a_2(s-3)z\big)-
$$
$$
B_s\big(\frac{1}{6}a_2(s+1)(6a_0a_2(4g(1+s)-8+4g^2(1+s)-s(2+3(s-3)s))+
$$
$$
a_2^2(-(s-2)(-6+(s-3)(s-1)^2s)+2g(6+s(9+4(s-4)s))+
$$
$$
2g^2(6+s(9+4(s-4)s)))+24(s-2)z)\big)+
$$
$$
\sum_{m=1}^{2g+4}\bigg((-1)^m\big(C_{s+m}^m(2a_0+5a_2)(z-a_0^2)+C_{s+m}^{m+1}6a_2(a_0^2-z)+C_{s+m}^{m+2}a_2(5a_2^2g(1+g)+
$$
$$
2a_0a_2(-5+g+g^2)-6a_0^2+2z)+C_{s+m}^{m+3}6a_2^2(2a_0-a_2g(1+g))+
$$
$$
C_{s+m}^{m+4}a_2^2(a_2(-5+2g(1+g))-6a_0)+C_{s+m}^{m+5}6a_2^3-C_{s+m}^{m+6}2a_2^3\big)-
$$
$$
\big(C_{s+m}^{m}(2a_0+a_2)(a_2^2(g^2+g+4)+3a_0^2+10a_0a_2+z)+C_{s+m}^{m+1}2a_2(9a_0^2+
$$
$$
2a_2^2+2a_0a_2(g^2+g+4)-z)+C_{s+m}^{m+2}a_2(18a_0^2-a_2^2(g^2+g-2)+2a_0a_2(g^2+g+19)+2z)+
$$
$$
C_{s+m}^{m+3}2a_2^2(18a_0+a_2g(1+g))-C_{s+m}^{m+4}a_2^2(18a_0+a_2(15+2g(1+g)))+
$$
$$
C_{s+m}^{m+5}18a_2^3+C_{s+m}^{m+6}6a_2^3\big)+2^m\big(C_{s+m}^{m}(2a_0+a_2)(a_0^2+2a_0a_2-a_2^2(g^2+g-1)-z)-
$$
$$
C_{s+m}^{m+1}4a_2(3a_0^2+a_2^2-2a_0a_2(g^2+g-2)+z)+C_{s+m}^{m+2}4a_2(6a_0^2+a_2^2g(g+1)-
$$
$$
2a_0a_2(g^2+g-5)-2z)-C_{s+m}^{m+3}16a_2^2(a_2g(g+1)-6a_0)+C_{s+m}^{m+4}16a_2^2(6a_0+
$$
$$
a_2(5-2g(g+1)))+C_{s+m}^{m+5}192a_2^3+C_{s+m}^{m+6}128a_2^3\big)\bigg)B_{s+m},
$$
where $0\leq s<2g+4, \quad {\rm C}_p^k=\frac{p!}{k!(p-k)}$ at $p\geq k,\ {\rm C}_p^k=0$ at $p<k.$ 

From the above formula it follows that
$$
\vartheta_{2g+8}=\vartheta_{2g+7}=\vartheta_{2g+6}=\vartheta_{2g+5}=\vartheta_{2g+4}=0,
$$
automatically. For example,
$$
\vartheta_{2g+8} = 12 a_2^3 (g+2) (4 g+5) B_{2 g+5} +8 a_2^3 (2 g+3) (2 g+5) (2 g+9) B_{2 g+6} +
$$
$$
\frac{2}{3} a_2 \big(6 a_0 a_2 (8 g^3+ 88 g^2+ 301 g+305) -a_2^2 (16 g^5+40 g^4 - 796 g^3- 4936 g^2 -
$$
$$
10275 g-7230) - 12 (2 g+5) z \big) B_{2 g+7} +\frac{2}{3} a_2 (2 g+9) (12 a_0 a_2 (2 g^3+26 g^2 + 
$$
$$
104 g +123) - a_2^2 (8 g^5 + 44 g^4 - 185 g^3- 1751 g^2 - 4017 g-2931) - 12 (g+3) z) B_{2 g+8}+
$$
$$
\sum_{m=1}^{2g+4}\bigg((-1)^m\big(C_{2g+8+m}^m(2a_0+5a_2)(z-a_0^2)+C_{2g+8+m}^{m+1}6a_2(a_0^2-z)+
$$
$$
C_{2g+8+m}^{m+2}a_2(5a_2^2g(g+1)+2a_0a_2(g^2+g-5)-6a_0^2+2z)+C_{2g+8+m}^{m+3}6a_2^2(2a_0-a_2g(g+1))+
$$
$$
C_{2g+8+m}^{m+4}(a_2^2(-6a_0+a_2(-5+2g(1+g))))+C_{2g+8+m}^{m+5}6a_2^3-C_{s+m}^{m+6}2a_2^3\big)-
$$
$$
\big(C_{2g+8+m}^{m}(2a_0+a_2)(a_2^2(g^2+g+4)+3a_0^2+10a_0a_2+z)+
$$
$$
C_{2g+8+m}^{m+1}2a_2(9a_0^2+2a_2^2+2a_0a_2(g^2+g+4)-z)+C_{2g+8+m}^{m+2}a_2(18a_0^2+
$$
$$
2a_0a_2(19+g+g^2)-a_2^2(g^2+g-2)+2z)+C_{2g+8+m}^{m+3}2a_2^2(18a_0+a_2g(g+1))-
$$
$$
C_{2g+8+m}^{m+4}a_2^2(18a_0+a_2(15+2g(g+1)))+C_{2g+8+m}^{m+5}18a_2^3+C_{2g+8+m}^{m+6}6a_2^3\big)+
$$
$$
2^m\big(C_{2g+8+m}^{m}(2a_0+a_2)(a_0^2+2a_0a_2-a_2^2(g^2+g-1)-z)-
$$
$$
C_{2g+8+m}^{m+1}4a_2(3a_0^2+a_2^2-2a_0a_2(g^2+g-2)+z)-C_{2g+8+m}^{m+2}4a_2(6a_0^2+a_2^2g(g+1)-
$$
$$
2a_0a_2(g^2+g-5)-2z)+C_{2g+8+m}^{m+3}16a_2^2(6a_0-a_2g(1+g))+C_{2g+8+m}^{m+4}(16a_2^2(6a_0+
$$
$$
a_2(5-2g(1+g))))+C_{2g+8+m}^{m+5}192a_2^3+C_{2g+8+m}^{m+6}128a_2^3\big)\bigg)B_{2g+8+m} = 0
$$
since we assume that $B_i=0$ if $i\geq 2g+3$. Similarly one can show that $\vartheta_{2g+7}=\vartheta_{2g+6}=\vartheta_{2g+5}=\vartheta_{2g+4}=0.$ So, we have
$$
R_n=\vartheta_{2g+3} n^{2g+3}+\vartheta_{2g+2} n^{2g+2}+\ldots+\vartheta_{0}=0.
$$
From the condition $\vartheta_{2g+3}=0$ we express $B_{2g}$ via $B_{2g+2}$
$$
B_{2 g} = B_{2 g+2} \frac{12 a_0 a_2 (2 g^2+2 g-1) - a_2^2 g (8 g^3+ 20 g^2 + 7 g-8) -
 12 z}{12 a_2^2 (2 g-1)}.
$$
We have
$$
R_n+R_{-n-1} = 2\vartheta_{2g+2} n^{2g+2}+(2g+2)\vartheta_{2g+2} n^{2g+1}+
$$
$$
(2\vartheta_{2g}-(2g+1)\vartheta_{2g+1}+(2g+2)(2g+1)\vartheta_{2g+2})n^{2g}+\ldots+(\vartheta_{2g+2}+\ldots+2 \vartheta_{0}) = 0.
$$
Hence from the condition $\vartheta_{2g+3}=0$ it follows that $\vartheta_{2g+2}=0$ automatically. Similarly, from $\vartheta_{2g-(2k-1)}=0$ we express $B_{2g-(2k+2)}$ via $B_{2g+2}$ and
$$
R_n+R_{-n-1} = 2\vartheta_{2g-2k} n^{2g-2k}+(2g-2k)\vartheta_{2g-2k} n^{2g-(2k+1)}+\ldots+(\vartheta_{2g+2k}+\ldots+2 \vartheta_{0}) = 0.
$$
Hence $\vartheta_{2g-2k}=0$ automatically.

So, from
$$
\vartheta_{2g+3}=\vartheta_{2g+1}=\ldots=\vartheta_{3}=\vartheta_{1}=0
$$
we express $ B_{2g}, \ldots, B_2, B_0$ via $B_{2g+2},$ herewith we have automatically
$$
\vartheta_{2g+2}=\vartheta_{2g}=\ldots=\vartheta_{2}=\vartheta_{0}=0.
$$
Next, we choose $B_{2g+2}$ such that coefficient at $z^g$ in $S_n(z)$ is $-U_n.$
We find $S_n,$ satisfying the equation (\ref{eq3}). Theorem 4 is proved.

\begin{rema}
We checked that for $g=1,\ldots,5$ the operator
$$
L_2=(T+\alpha_2n^2+\alpha_1n+\alpha_0)^2-g(g+1)\alpha_2n(\alpha_2n+\alpha_1), \ \alpha_2\neq0
$$
commute with $L_{2g+1}.$ Probably this operator $L_2$ commutes with $L_{2g+1}$ for all $g \in \mathbb{N}, \alpha_2 \neq 0, \alpha_1, \alpha_0.$
\end{rema}

\noindent{\bf Example 3}. At $g=1$ the operator
$$
L_2=(T + a_2 n^2 + a_0)^2 - 2 a_2^2 n^2
$$
commutes with the operator $L_3$
$$
L_3 = T^3 + (a_2 (3 n^2+6 n+5)+3 a_0) T^2 +
$$
$$
\frac{1}{4} (a_2 (2 n^2 + 2 n+1)+2 a_0 ) (a_2 (6 n^2 + 6 n-1) + 6 a_0) T +
$$
$$
\frac{1}{4} (a_2 (2 n^2 -2 n-1)+ 2 a_0) (a_2 (n^2-1)+a_0) (a_2 (2 n^2 +2 n -1)+2 a_0),
$$
$$
S_n=a_2^3 n^4+\frac{1}{4} a_2 (12 a_0 a_2 - 9 a_2^2 - 4 z)n^2+\frac{1}{4} (8 a_0^2 a_2 - 5 a_0 a_2^2 + a_2^3 - 4 a_0 z),
$$
$$
Q_n=z-\frac{1}{4} a_2 (8 a_0 + a_2 (4 n (n+1)-3 )).
$$
The spectral curve is given by the equation
$$
w^2 = \frac{1}{16} (z + a_2^2-2 a_0 a_2) (4 z +a_2^2-4 a_0 a_2)^2.
$$


\subsection{Theorem 5}
Let
$$
U_n = \beta a^n, \qquad W_n = \frac{a^{2 g} + a^{2 g+2} - a^{4 g+2} - 1}{(a^{2 g+1} + 1)^2} \beta^2 a^{2 n}.
$$
The equation (\ref{eq3}) has the form
$$
S_{n - 1} ((a + 1) (a^{2 g + 2} + a)^2 \beta z-a^{2 (n + g + 1)} (a + 1)^3 \beta^3)+
$$
$$
S_n a (a + 1) \beta (a^{2 n} (a + 1)^2 (a^{2 g + 1} + a^{4 g + 2} + 1) \beta^2 +a (a^{2 g + 1} + 1)^2 z)-
$$
$$
S_{n + 1} (a^{2 n + 1} (a + 1)^3 (a^{2 g + 1} + a^{4 g + 2} + 1) \beta^3 + (a + 1) (a^{2 g + 1} + 1)^2 \beta z)
$$
\begin{equation}
\label{eq11}
S_{n + 2} (a^{2 (n + g + 1)} (a + 1)^3 \beta^3 - (a + 1) (a^{2 g + 1} + 1)^2 \beta z) = 0.
\end{equation}
Let
$$
S_n=G_{2g+1}a^{(2g+1)n}+G_{2g-1}a^{(2g-1)n}+\ldots+G_1a^n, \qquad G_i=G_i(z).
$$
After substituting $S_n$ in (\ref{eq11}) we obtain
$$
R_n = \mu_{2g+4} a^{(2g+4)n}+\mu_{2g+2} a^{(2g+2)n}+\ldots+\mu_{4} a^{4n} = 0,
$$
where
$$
\mu_p =  G_{p-3} \frac{(a+1)^2 (a^3 - a^p) (a^{2 g+4} - a^p) (a^{2 g + p}-a^2) \beta^2}{(a^{2 g+1}+1)^2} +
$$
$$
G_{p-1} a^2 (a^2 - a^p) (a + a^p) (a^2 + a^p) z .
$$
From the above formula it follows that
$$
\mu_{2g+4} = G_{2 g+3} a^2 (a^2 - a^{2g+4}) (a + a^{2g+4}) (a^2 + a^{2g+4}) z
$$
(we assume that $G_{2g+3}=0$). So, we have $g$ equations $\mu_{2g+2}=\mu_{2g}=\ldots=\mu_{4}=0$ and $g+1$ unknown functions $G_1, G_3,\ldots, G_{2g+1}.$ From $\mu_p=0,$
we express $G_{2g-1}, G_{2g-3},\ldots, G_1$ via $G_{2g+1}.$ For example,
$$
G_{2g-1} =-G_{2 g+1}\frac{a^{1 - 2 g} (a^{2 g+1}+1)^3 z}{\beta_1^2(a-1) (a+1)^3 (a^{2 g}-a)}.
$$
Next, we choose $G_{2g+1}$ such that coefficient at $z^g$ in $S_n(z)$ is $-U_n.$
We find $S_n,$ satisfying the equation (\ref{eq3}). Theorem 5 is proved.

\noindent{\bf Example 4}.
At $g=1$ the operator
$$
L_2=(T+\beta a^n)^2-\frac{(a^2 + a^4 - a^6-1)\beta^2}{(a^3+1)^2} a^{2 n}
$$
commutes with the operator $L_3$
$$
L_3 = T^3+(a^2 +a +1) \beta a^n T^2+\frac{(a^2+a+1) \beta^2 a^{2 n+1}}{a^2-a+1}T+\frac{\beta^3 a^{3 n+3}}{(a^2-a+1)^3}.
$$
$$
S_n=- \beta a^n z+\frac{(a-1)^2 \beta^3 a^{3 n+2}}{(a^2 - a+1)^3},
\qquad
Q_n = z - \frac{(a-1)^2 \beta^2a^{2 n}}{(a^2-a+1)^2}.
$$
The spectral curve is given by the equation
$$
w^2= z^3.
$$

\section{One--point higher rank positive operators}

In this section we introduce spectral data for positive operators of rank $l>1.$

We take the following spectral data
$$
S=\{\Gamma, q, k^{-1}, K_n, \xi_0(n), \gamma, \alpha \},
$$
where $\Gamma$ is a Riemann surface of genus $g,$ $q \in \Gamma$ is a fixed point, $k^{-1}$ is a local parameter near $q,$ $K_n \in \Gamma, n \in \mathbb{Z}$ is a set of points in general position, $\xi_0(n)=(\xi^1_0(n),\ldots,\xi^{l-1}_0(n),1),$ at $n>0,$ $\xi_0(n)=(1,\xi^2_0(n),\ldots,\xi^{l}_0(n)),$ at $n<0,$ $\xi^j_0(n)$ is a function in $n,$ for simplicity we assume $\xi^j_0(n)\neq 0,$ $\gamma = \gamma_1+\ldots+\gamma_{lg}$ is a divisor ($\gamma_j \in \Gamma$  are in general position), $\alpha$ is a set of vectors
$$
\alpha_1,\ldots,\alpha_{lg}, \qquad \alpha_j = (\alpha_{j,1},\ldots,\alpha_{j,l-1}).
$$
Pair $(\gamma, \alpha)$ is called {\it Tyurin parameters}, $(\gamma, \alpha)$ defines a semistable holomorphic bundle over $\Gamma$ with holomorphic sections $\eta_1,\ldots,\eta_l.$ Points $\gamma_1,\ldots,\gamma_{lg}$ are points of linear dependence of the sections, herewith
$$
\eta_l(\gamma_k) = \sum^{l-1}_{i=1}\alpha_{j,i}\eta_j(\gamma_k), \qquad k=1,\ldots,lg.
$$
\begin{theorem}
There is a unique vector function $\psi(n,P)=(\psi_1(n,P),\ldots,\psi_l(n,P))$ which satisfies the following properties.

\noindent{\bf 1}. In the neighborhood of point $q$ the function $\psi(n,P)$ has the form
\begin{equation}
\label{eq12}
\psi(n,P) = \bigg(\sum_{s=0}^{\infty}\frac{\xi_s(n)}{k^s} \bigg)
\begin{pmatrix}
0 & 1 & 0 & \ldots & 0 & 0 \\
0 & 0 & 1 & \ldots & 0 & 0 \\
\hdotsfor{6} \\
0 & 0 & 0 & \ldots & 0 & 1 \\
k & 0 & 0 & \ldots  & 0 & 0 \\
\end{pmatrix}^n,
\end{equation}
$$
\xi_s(n) = (\xi^1_s(n),\ldots,\xi^{l}_s(n)).
$$

\noindent{\bf 2}. The vector function $\psi(n,P)$ has $gl$ simple poles in $\gamma_1,\ldots,\gamma_{lg}$ such that
\begin{equation}
\label{eq13}
Res_{\gamma_i} \psi_j = \alpha_{i,j} Res_{\gamma_i} \psi_l.
\end{equation}

\noindent{\bf 3}. Let $n>0$ and let $n=l m+s, \ m \in \mathbb{N}, \ 0<s\leq l,$ $\psi$ has simple zeros in $K_1,\ldots,K_m$
\begin{equation}
\label{eq14}
\psi(K_p) = 0, \qquad 1\leq p \leq m,
\end{equation}
additionally $\psi_1,\ldots,\psi_s$ have simple zero in $K_{m+1}$
\begin{equation}
\label{eq15}
\psi_1(K_{m+1})=0,\ldots,\psi_s(K_{m+1})=0.
\end{equation}
Let $n<0$ and let $n=- l m-s , \ m \in \mathbb{N}, \ 0<s\leq l,$ $\psi$ has simple poles in $K_{-m},\ldots,K_{-1}$
\begin{equation}
\label{eq16}
\psi(K_p)=\infty, \qquad -m\leq p \leq -1,
\end{equation}
additionally $\psi_1,\ldots,\psi_s$ have simple pole in $K_{-m-1}$
\begin{equation}
\label{eq17}
\psi_1(K_{-m-1})=\infty,\ldots,\psi_s(K_{-m-1})=\infty.
\end{equation}

\noindent{\bf 4}. $\psi(0,P) = \xi_0(0).$
\end{theorem}
\noindent\textbf{Proof.} Let us consider the case $n>0.$
Let
$$
A=\begin{pmatrix}
0 & 1 & 0 & \ldots & 0 & 0 \\
0 & 0 & 1 & \ldots & 0 & 0 \\
\hdotsfor{6} \\
0 & 0 & 0 & \ldots & 0 & 1 \\
k & 0 & 0 & \ldots  & 0 & 0 \\
\end{pmatrix}.
$$
Then
$$
A^2=\begin{pmatrix}
0 & 0 & 1 & 0 & \ldots & 0 & 0 \\
0 & 0 & 0 & 1& \ldots & 0 & 0 \\
\hdotsfor{7} \\
0 & 0 & 0 & 0 & \ldots & 0 & 1 \\
k & 0 & 0 & 0 & \ldots & 0 & 0 \\
0 & k & 0 & 0 & \ldots & 0 & 0 \\
\end{pmatrix},
\quad
A^3=\begin{pmatrix}
0 & 0 & 0 & 1 & \ldots & 0 & 0 \\
0 & 0 & 0 & 0 & \ldots & 0 & 0 \\
\hdotsfor{7} \\
k & 0 & 0 & 0 & \ldots & 0 & 0 \\
0 & k & 0 & 0 & \ldots & 0 & 0 \\
0 & 0 & k & 0 & \ldots & 0 & 0 \\
\end{pmatrix}, \quad \ldots,
$$
$$
A^l=\begin{pmatrix}
k & 0 & 0 & \ldots & 0 & 0 \\
0 & k & 0 & \ldots & 0 & 0 \\
\hdotsfor{6} \\
0 & 0 & 0 & \ldots & k & 0 \\
0 & 0 & 0 & \ldots  & 0 & k \\
\end{pmatrix}.
$$
We have $A^n = k^m A^s.$
Hence, the pole divisors of $\psi_j$ are
$$
(\psi_j)_\infty = \gamma_1+\ldots+\gamma_{lg}+q(m+1), \qquad 0<j\leq s
$$
$$
(\psi_j)_\infty = \gamma_1+\ldots+\gamma_{lg}+qm, \qquad s<j\leq l.
$$
The conditions (\ref{eq12})--(\ref{eq15}) define the unique function $\psi(n,P).$

Similarly, when $n<0,$ the pole divisors of $\psi_j$ are
$$
(\psi_j)_\infty = \gamma_1+\ldots+\gamma_{lg}+K_{-1}+\ldots+K_{-m-1}, \qquad 0<j\leq s
$$
$$
(\psi_j)_\infty = \gamma_1+\ldots+\gamma_{lg}+K_{-1}+\ldots+K_{-m}, \qquad s<j\leq l.
$$
The conditions (\ref{eq12})--(\ref{eq13}), (\ref{eq16})--(\ref{eq17}) define the unique function $\psi(n,P).$
Theorem 7 is proved.

For meromorphic functions $f(P), g(P)$ with unique poles in $q$ of orders $m, s$ there are two commuting operators $L_{lm}, L_{ns}$ of orders $lm$ and $ns$ such that
$$
L_{lm} \psi = f(P) \psi, \qquad L_{ns} \psi = g(P) \psi.
$$
Theorem 7 is an analogue of the Krichever's Theorem (see \cite{K}) for spectral data of higher rank commuting ordinary differential operators. Higher rank commuting ordinary differential operators were studied for example in \cite{KN1}--\cite{Mir1} (see another citations in \cite{Mir}).

We plan to study commuting higher rank positive operators later. Here we give only one example of such operators.

\noindent{\bf Example 5}.
Operator
$$
L_4 = T^4+(a_2 n^2+a_1 n+a_0)T^3+(\frac{3}{8} (a_1 + a_2 (n-1)) n (2 a_0 + a_1 (n-1) +
$$
$$
a_2 (n^2-n-2)))T^2 - \frac{1}{16} (a_0 + a_1 (n-1) + a_2 (n-2) n) (2 a_0^2 - a_1^2 (n-2) n -
$$
$$
a_0 (a_1 + 2 a_2 (n-1)^2 + 2 a_1 n) - a_1 a_2 (2 n^3- 6 n^2 - n +2) - a_2^2 n (n^3 - 4 n^2 - n + 10)))T + $$
$$
\frac{1}{256} (a_1 + a_2 (n-3)) n (2 a_0 - 4 a_2 + (n-3) (a_1 + a_2 n)) (-4 a_0^2 + a_1^2 (n-2) (n-1) +
$$
$$
a_2^2 (n-2) (n-1) ((n-3) n-6) + 2 a_0 (a_1 n + a_2 ((n-3) n+4)) + a_1 a_2 (6 + n (5 + n (2 n-9)))),
$$
commutes with an operator of order 6, the spectral curve is given by the equation
$$
w^2 = \frac{1}{262144} \big(32 z + (a_0 - a_1) (a_0 - a_2) (a_0 (2 a_0 - a_1) - 2 (a_0 + a_1) a_2) \big)^2\times
$$
$$
\big(256 z + (-2 a_0^2 + a_1^2 + 4 a_0 a_2 + 3 (a_1 - 2 a_2) a_2)^2 \big).
$$
At $a_2=2, \ a_1=a_0=0$
$$
L_6 =T^6+(3 n^2 + 6 n +8) T^5+\frac{1}{4} (n (n+1) (32 + 15 n (n+1))-6)T^4+
$$
$$
\frac{1}{2} n^2 (n^2-2) (5 n^2+7)T^3+\frac{1}{16} (n-2) (n-1) n (n+1) ((n-1) n (15 (n-1) n -38)-36)T^2+
$$
$$
\frac{1}{16} (n-2)^2 n^2 (12 + (n-2) n ((n-2) n-5) (3 (n-2) n-11))T+
$$
$$
\frac{1}{64} (n-4) (n-3) (n-2) (n-1) n (n+1) ((n-3) n-6) ((n-4) (n-3) n (n+1)-6).
$$


\begin{thebibliography}{1}

\bibitem{K3} I.M.~Krichever, \textit{Algebraic curves and non--linear difference equations}, Russian Math. Surveys, {\bf 33}:4, (1978) 215--216.

\bibitem{Mum} D.~Mumford, {\it An algebro--geometric construction of commuting operators and of solutions to the Toda lattice equation, Korteweg--de Vries equation and related non--linear equations}, Proceedings of the International Symposium on Algebraic Geometry (Kyoto Univ., Kyoto, 1977), Kinokuniya, Tokyo, 1978, 115--153.

\bibitem{KN2} I.M.~Krichever, S.P.~Novikov, \textit{Two--dimensionalized Toda lattice, commuting difference operators, and holomorphic bundles},
Russian Math. Surveys, {\bf 58}:3, (2003) 473--510.

\bibitem{KN3} I.M.~Krichever, S.P.~Novikov, \textit{Holomorphic bundles and scalar difference operators: one-point constructions},
Russian Math. Surveys, {\bf 55}:1, (2000) 180--181.


\bibitem{MM1}
G.S.~Mauleshova, A.E.~Mironov, \textit{Commuting difference operators of rank two},
Russian Math. Surveys, {\bf 70}:3, (2015) 557--559.

\bibitem{MM2} G.S.~Mauleshova, A.E.~Mironov, \textit{One--point commuting difference operators of rank 1}, Doklady Mathematics, {\bf 93}:1, (2016) 62--64.

\bibitem{K5}
I.M.~Krichever, \textit{Two--dimensional periodic difference operators and algebraic geometry},
Dokl. Akad. Nauk SSSR, {\bf 285}:1, (1985) 31--36.

\bibitem{MM3}
G.S.~Mauleshova, A.E.~Mironov, \textit{One--point commuting difference operators of rank 1 and their relation with finite--gap Schr\"{o}dinger operators}, Doklady Mathematics, {\bf 97}:1, (2018) 62--64.

\bibitem{K}
I.M.~Krichever, \textit{Commutative rings of ordinary linear differential operators}, Funct. Anal. Appl., {\bf 12}:3, (1978) 175–185.


\bibitem{KN1}
I.M.~Krichever, S.P.~Novikov, \textit{Holomorphic bundles over algebraic curves and non--linear equations}, Russian Mathematical Surveys, {\bf 35}:6, (1980) 53--79.

\bibitem{G}
P.G~Grinevich, \textit{Rational solutions for the equation of commutation of differential operators}, Functional Analysis and Its Applications, {\bf 16}:1, (1982) 15--19

\bibitem{MokhO}
O.I.~Mokhov, \textit{Commuting differential operators of rank 3, and nonlinear differential equations}, Mathematics of the USSR--Izvestiya, {\bf 35}:3, (1990) 629--655

\bibitem{PW}
E.~Previato, G.~Wilson, \textit{Differential operators and rank $2$ bundles over elliptic curves}, Compositio Math., {\bf 81}:1, (1992) 107--119.

\bibitem{LP}
G.~Latham, E.~Previato,  \textit{Darboux transformations for higher-rank Kadomtsev--Petviashvili and Krichever--Novikov equations},
Acta Appl. Math., {\bf 39}, (1995) 405--433.

\bibitem{Mir1}
A.E.~Mironov, \textit{Self--adjoint commuting differential operators.}, Inventiones mathematicae, {\bf 197}:2, (2014) 417--431.

\bibitem{Mir}
A.E.~Mironov, \textit{Self--adjoint commuting differential operators of rank two},
Russian Mathematical Surveys, {\bf 71}:4, (2016) 751--779.



\end{thebibliography}
\end{document}